\documentclass{amsart}
\usepackage{amsthm}%
\usepackage{hyperref}
\usepackage{amsrefs}
\usepackage{comment}
\usepackage{float}
\usepackage{caption}

\usepackage{graphicx}%
\usepackage{listings}%
\usepackage{seqsplit}
\newtheorem{theorem}{Theorem}%  meant for continuous numbers
\usepackage{hyperref}
\usepackage{comment}
\usepackage{graphicx} %Inserting graphics
\usepackage{float} %Forcing exact graphic placement - H

\newcommand{\etal}{\textit{et~al.\!}}
\newcommand{\m}{{\mathbf m}}
\newcommand{\n}{{\mathbf n}}

\newcommand{\T}{{\mathbf t}}
\newcommand{\Sb}{{\mathbf S}}

\newcommand{\Gb}{{\mathbf G}}

\newcommand{\mC}{C_{\m}}

\newcommand{\BigComment}[1]{}

\providecommand{\U}[1]{\protect\rule{.1in}{.1in}}

\newcommand{\WR}{Wildberger and Rubine}

\newcommand{\code}[1]{\texttt{#1}}

\begin{document} %----------------------------------------------------------------------------

\title{Computing the 4D Geode}
\markright{Computing the 4D Geode}
\author{Dean Rubine}
\address{Lee, New Hampshire}
\email{DeanRubineMath@gmail.com}
\date{\today}

\begin{abstract}
The closed form for the hyper-Catalan number  $C[m_2,m_3,m_4,\ldots]$,
which
counts the number of subdivisions of a roofed
polygon into $m_2$ triangles, $m_3$ quadrilaterals,  $m_4$ pentagons,
etc., has been known since % Erd\'elyi and Etherington, 
1940.
In 2025,
\WR{} showed
its generating series $\Sb[t_2,t_3,t_4,\ldots]$ is a zero of the general geometric univariate polynomial.
They note the factorization $\Sb-1=(t_2 + t_3 + t_4 + \ldots)\Gb$, where the factor $\Gb$ is called the Geode.
Later in 2025, Amderberhan, Kauers and Zeilberger issued a challenge to compute $G[1000,1000,1000,1000]$, the coefficient of $t_2^{1000}t_3^{1000}t_4^{1000}t_5^{1000}$ in $\Gb$.
The reward is a donation to OEIS.
We describe the computation, give the value and claim the reward.
\end{abstract}

% 
% 05A15     Exact enumeration problems, generating functions [See also 33Cxx, 33Dxx]
% 08-02  	Research exposition (monographs, survey articles) pertaining to general algebraic systems
% 17-02  	Research exposition (monographs, survey articles) pertaining to nonassociative rings and algebras
% 12D10  	Polynomials in real and complex fields: location of zeros (algebraic theorems) {For the analytic theory, see 26C10, 30C15}
% 40-02  	Research exposition (monographs, survey articles) pertaining to sequences, series, summability

% KEYWORDS:  polynomial solutions, polynomial roots, power series zeros, hyper-Catalan numbers, subdigon polyseries, Geode array
\date{\today}
%author{\fnm{Dean} \sur{title}} \email{DeanRubineMath@gmail.com, orcid.org/0009-0000-8767-8137}
\maketitle
%\pacs[JEL Classification]{D8, H51}
%\pacs[MSC Classification]{05A15}

\section{Introduction}

As a computer scientist and coauthor of \textit{A Hyper-Catalan Series Solution to Polynomial Equations, and the Geode}~\cite{Wildberger2025},
I was in a very good position to attempt the \textbf{Geode Challenge}, a computation problem posed by Amderberhan, Kauers and Zeilberger in August, 2025~\cite{Amdeberhan2025}.
I had already written a paper~\cite{rubine2025hyper} that developed a recurrence (Theorem 8 therein) that could straightforwardly calculate Geode values for arbitrary type vectors.
That program (figure \ref{fig:code1} below) was not particularly efficient, but it finished quickly enough on the simple type vectors I was using for pedagogical examples.

The inefficiency never mattered until the Geode Challenge.
I resisted trying for a while, perhaps because I was not that familiar with the mathematics used in the two Amderberhan \etal{} papers~\cites{Amdeberhan2025,Zeilberger2025}.
Then on a weekend in mid-November, it occurred to me that I could speed up my simple Geode calculator with a simple cache of some computational results, which led to me successfully computing $G[200,200,200,200]$ and $G[250,250,250,250]$, a quarter of the way to the prize.
I sent my results to the authors of the challenge, who were very encouraging.
After another couple of rounds of speedups, I finished computing $G[1000,1000,1000,1000]$ on the solstice, December 21, 2025.
This paper tells the story.

\section{The Geometric Polynomial Formula}

\begin{figure}
    \centering
    \includegraphics[width=1\linewidth]{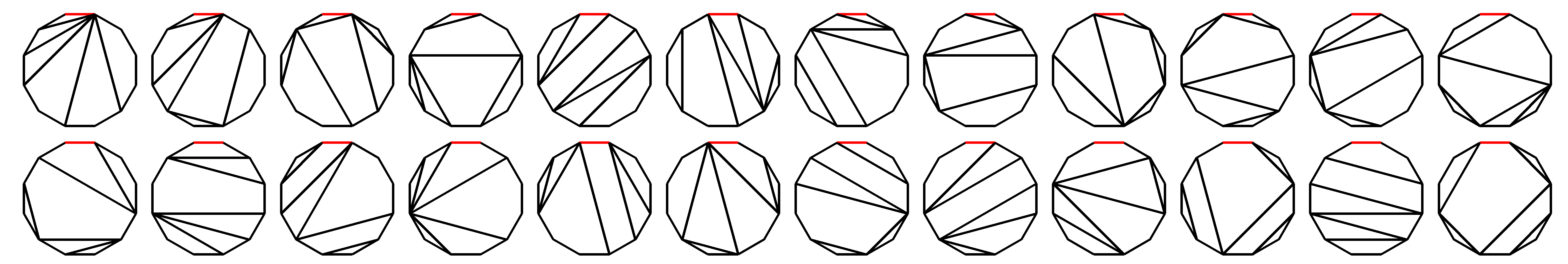}
    \caption{A few of the $C[3,2,1]=43680$ subdigons subdivided into three triangles, two quadrilaterals and a pentagon.}
    \label{fig:C321}
\end{figure}

Erd\'elyi and Etherington, 1940~\cite{Erdelyi1940} appear to be the first to give a combinatorial interpretation of the hyper-Catalan coefficient $\mC=C[m_2, m_3, m_4,\ldots] $.
The combinatorial object being counted is a roofed polygon subdivided by non-crossing diagonals, which we call a \textbf{subdigon}.
$\mC$ counts the number of subdigons with $m_2$ triangular faces, $m_3$ quadrilateral faces, etc., i.e. $\mC$ is the number of subdigons of \textbf{type} $\m=[m_2, m_3, m_4, \ldots]$.
Appending zeros does not change the type.
Figure \ref{fig:C321} illustrates some subdigons of type $[3,2,1]$. 

The indexing starting at two is a bit awkward; Amderberhan \textit{et al.}~\cite{Amdeberhan2025} renotate to start at one.  Here we retain the indexing starting at two, as there are times in which we need the array slightly more general than the hyper-Catalans (see for example
\cites{Rubine2025Finite,Rubine2025Fine}).
That array includes $m_1$, which somewhat problematically refers to 2-gons in the context of subdivided polygons.
The generalized array also counts plane trees of a type; there $m_1$ refers to the number of unary nodes.

Our computation of $\Gb[1000,1000,1000,1000]$ requires computing over a billion hyper-Catalan numbers, so the details of their calculation will be important to efficiency.
As our concern is the four dimensional Geode, we may restrict ourselves to fifth degree polynomials, which means we are only concerned with non-zero $m_2, m_3, m_4,$ and $ m_5$, subdivisions into triangles, quadrilaterals, pentagons and hexagons only.  For the rest of this paper $\m$ will always be a vector of four elements, $\m=[m_2, m_3, m_4, m_5]$.  We will restrict our theorems to this special case.

\begin{theorem}[4D Hyper-Catalan Closed Form, Erd\'elyi and Etherington, 1940%~\cite{Erdelyi1940}
] \label{thm:hypercatalanclosedform}
The number of subdigons of type $\m=[m_2, m_3, m_4, m_5]$ is
\begin{align*}
\mC &=\dfrac{( 2m_2 + 3m_3 + 4m_4 + 5m_5)!}{(1 + m_2 + 2m_3 + 3m_4 + 5m_5)! \,  m_2! \, m_3! \,  m_4! \, m_5!} 
=  \dfrac{(E_\m - 1)!}{(V_\m - 1)! \ \m!}
.
\end{align*}
where $E_\m \equiv 1+2m_2 + 3m_3 + 4m_4 + 5m_5$ is the number of edges in a subdigon of type $\m$, $V_\m \equiv 2+m_2 + 2m_3 + 3m_4 + 5m_5$ is the number of vertices in a subdigon of type $\m$, and $\m! \equiv m_2! \, m_3! \,  m_4! \, m_5!$.
\end{theorem}

Wildberger and Rubine's~\cite{Wildberger2025} main theorem is that the generating series for the hyper-Catalan numbers is a formal power series solution of the \textbf{general geometric polynomial}, which means general except for a constant of $1$ and a linear coefficient of $-1$.
\begin{theorem}[The geometric quintic formula, \WR{}, 2025] \label{thm:geompoly}
The generating sum for the four dimensional hyper-Catalan numbers:
\begin{align}
\Sb \equiv \sum_{m_2 \ge 0} \sum_{m_3 \ge 0} \sum_{m_4 \ge 0} \sum_{m_5 \ge 0} 
% \Sb \equiv \sum_{m_2,m_3,m_4, m_5 \ge 0} 
C[m_2, m_3, m_4,m_5] t_2^{m_2} t_3^{m_3} t_4^{m_4} t_5^{m_5} \equiv \sum_{\m \ge 0} \mC \T^{\m}
\end{align}
is a formal series zero of the general geometric quintic polynomial:
\begin{align} \label{eqn:galpha}
g(\alpha) = 1 - \alpha + t_2 \alpha^2 + t_3 \alpha^3 + t_4 \alpha^4 + t_5 \alpha^5
\end{align}
\end{theorem}
Please see their paper for the generalization, detailed definitions and the proof.

\section{The Geode}
\WR{} uncover the Geode while examining a face layering of $\Sb$, which they call $ \Sb_F  = \Sb[ft_2, ft_3, ft_4, ft_5]$ (here restricted to four dimensions).
Expanding, they get a layering of $\Sb$ into ongoing power series
of a given total degree, accounting for subdigons with a given number of faces, and they notice a factorization at each face level:
\begin{align}
\nonumber \Sb_{F}  = {} & 1+(  t_{2}+t_{3}+t_{4}+t_5)  f
\\[-2pt] \nonumber &
+ (2 t_2^2 + 5 t_2 t_3 + 6 t_2 t_4 + 7 t_2 t_5 + 3 t_3^2 + 7 t_3 t_4 + 8 t_3 t_5 + 4 t_4^2 + 9 t_4 t_5 + 5 t_5^2) f^2 
\\[-2pt] \nonumber &
+ (5 t_2^3 + 21 t_2^2 t_3 + 28 t_2^2 t_4 + 36 t_2^2 t_5 + 28 t_2 t_3^2 + 72 t_2 t_3 t_4 + 90 t_2 t_3 t_5 + 45 t_2 t_4^2 
\\[-4pt] \nonumber & \qquad
+ 110 t_2 t_4 t_5 + 66 t_2 t_5^2 + 12 t_3^3 + 45 t_3^2 t_4 + 55 t_3^2 t_5 + 55 t_3 t_4^2 + 132 t_3 t_4 t_5 
\\[-4pt]  & \qquad
+ 78 t_3 t_5^2 + 22 t_4^3 + 78 t_4^2 t_5 + 91 t_4 t_5^2 + 35 t_5^3) f^3 + \ldots
\\ \nonumber
\\ \nonumber
\Sb_{F}  =  1  & +  \Sb_1  f + \Sb_1  (2 t_2 + 3 t_3 + 4 t_4 + 5 t_5) f^2
\\[-2pt] \nonumber &
 +\Sb_1  (5 t_2^2 + 16 t_2 t_3 + 23 t_2 t_4 + 31 t_2 t_5 + 12 t_3^2 + 33 t_3 t_4 + 43 t_3 t_5
\\[-4pt]  & \qquad
+ 22 t_4^2 + 56 t_4 t_5 + 35 t_5^2) f^{3}
+\cdots 
\end{align}
where $\Sb_1 \equiv t_2 + t_3 + t_4 + t_5$.
That leads them to:

\begin{theorem}[The Geode factorization, \WR{}, 2025] \label{thm:geodefactorization}
There is a unique polyseries $\Gb $ satisfying $\Sb-1=\Sb_1\Gb.$
\end{theorem}
That concludes the summary of Wildberger and Rubine's 2025 paper.
We refer readers there for the inductive proof of the Geode factorizaton.

Much of the research attention since the publication of Wildberger and Rubine~\cite{Wildberger2025} has been focused on the Geode~\cites{Amdeberhan2025,rubine2025hyper, rubine2025exer, Gessel2025, Zeilberger2025, Gossow2025}.
In particular, the conjectures in the paper have all been resolved (all but one true), and combinatorial interpretations of the Geode have been given by Gessel~\cite{Gessel2025} in terms of lattice paths, and by Gossow~\cite{Gossow2025} in terms of subdivided polygons. 

The closed form for the general Geode element remains an open problem.
The current author believes that such a form is unlikely to exist.

\section{An OEIS Tangent}

The challenge from Amderberhan, Kauers and Zeilberger~\cite{Amdeberhan2025} is to compute $G[1000,1000,1000,1000] = [ t_2^{1000}  t_3^{1000} t_4^{1000} t_5^{1000} ] \Gb$.

A straightforward calculation of $G[1000,1000,1000,1000]$ from its \mbox{definition} involves generating the sum of all the terms of $\Sb[t_2, t_3, t_4, t_5]$ up to say $C[1001,1001, $ $ 1001, 1001] t_2^{1001} t_3^{1001}  t_4^{1001} t_5^{1001}$, dividing that (without the constant) by $t_2+t_3+t_4+t_5$, and picking out the coefficient of the appropriate term.
That is \textit{prima facie} an onerous calculation, involving the construction and division of a four variable polynomial with over a trillion terms.

Actually, we saw above that $\Sb_F$ factors at each face level, where the number of faces in a subdigon of type $\m$ is $F_\m=m_2+m_3+m_4+m_5$, the degree of the monomial $t_2^{m_2} t_3^{m_3} t_4^{m_4} t_5^{m_5}$.
So we can reduce our computation task by focusing on the particular face level of interest, here $4001$ faces (because $\Sb_1$ is of degree $1$,
and the term we seek in the Geode factorization comes from particular hyper-Catalan terms one degree higher.
The terms of $\Sb[t_2,t_3,t_4,t_5]$ corresponding to subdigons with $4001$ faces (i.e. of degree $4001$) is still a very large polynomial.

How large?  How many terms are in each face level of $\Sb[t_2,t_3,t_4,t_5]$?
The reward for the Geode challenge is a donation to OEIS~\cite{OEIS} in one's honor, so let's take a little OEIS tangent and see if we can figure this out.

We are trying to count how many monomials of the form $t_2^{m_2} t_3^{m_3} t_4^{m_4} t_5^{m_5}$  there are of degree 4001. Those are the 4-tuples $[m_2,m_3,m_4,m_5]$ where $F_\m=m_2+m_3+m_4+m_5=4001$.   
A short program $a(n)=\sum_{F_\m=n} 1\ $ generates 1, 4, 10, 20, 35, 56, 84, 120, 165, 220, 286, 364, 455, 560, 680, 816, 969, 1140, 1330, 1540, ... which OEIS happily informs us is A292, the tetrahedral numbers, a binomial coefficient, $\binom{n+3}{3}$ for our 0-indexed list. So to directly compute $G[1000,1000,1000,1000]$ we need generate a face layer of $\Sb$ with $\binom{4001+3}{3}=10,690,684,004$ terms, which we then divide by $t_2 + t_3 + t_4 + t_5$. Thanks, OEIS!

Ten billion terms might not be out of the question; let's call that future work.
Below, we take a different path.

A related question is how many subdigons with a four element type vector have 4001 faces?
The generating series for that is
\begin{align} \nonumber 
\Sb[f,f,f,f] &
= 4 f +56 f^{2} + 1104 f^{3} +  25376 f^{4} + 636480 f^{5} +16890240 f^{6} 
\\[-2pt] & 
+ 466254080 f^{7}  + 13251822080 f^{8} + 385188955136 f^{9} + % 11397332301824 f^{10} + 
\cdots 
\end{align}
No joy in OEIS for that one; we can try what is left after we factor out the $4f$.
\begin{align}
\Sb[f,f,f,f] & = 4f( 1 + 14 f +  276 f^{2} +6344 f^{3} +159120 f^{4}  + \cdots )
% 4 f (2849333075456 f^{9} + 96297238784 f^{8} + 3312955520 f^{7} + 116563520 f^{6} + 4222560 f^{5}
\end{align}
Still not in OEIS; so sad.

This particular sequence, from $\Sb[f,f,f,f]$, is of course dependent on our choice of four dimensions, i.e. of limiting ourselves to solving quintic equations. Let's see if the analogous sequences for fewer dimensions are in OEIS.  Yes, one dimension, degree $2$, is $\Sb[f]$, the generating series for the Catalan numbers, our old friend A108 in OEIS.
$\Sb[f,f]$ generates A27307, the 3-Schroeder numbers.
$\Sb[f,f,f]$ generates A371657.
I have not yet looked through the entries to see if there is a formula, and the degree five version we need is not in OEIS yet.
In this case the OEIS has failed to tell us how many subdigons divided into just 3-gons, 4-gons, 5-gons, and 6-gons have exactly 4001 faces.
But as usual, it has given us a few leads to follow.
(By the way, a program directly calculating $\sum_{F_\m=4001} \mC  $ never finished either.)

\begin{comment}

Now that we are thinking about it, we see it is only the four hyper-Catalan terms with three indices that are $1000$ and one index that is $1001$ that contribute:
\begin{align} \nonumber
(t_2+t_3 & +t_4+t_5) ( \ldots + G[1000,1000,1000,1000] t_2^{1000} t_3^{1000}  t_4^{1000} t_5^{1000} + \ldots)  = 
\\&
( \ldots   \nonumber
+ C[1001,1000,1000,1000]t_2^{1001} t_3^{1000}  t_4^{1000} t_5^{1000} 
\\& \ \ \ \ \ \; \nonumber
+ C[1000,1001,1000,1000]t_2^{1000} t_3^{1001}  t_4^{1000} t_5^{1000}
\\& \ \ \ \ \ \; \nonumber
+ C[1000,1000,1001,1000]t_2^{1000} t_3^{1000}  t_4^{1001} t_5^{1000}
\\& \ \ \ \ \ \;
+ C[1000,1000,1000,1001]t_2^{1000} t_3^{1000}  t_4^{1000} t_5^{1001}
+ \ldots)
\end{align}
\end{comment}

\section{The Geode Recurrence}

Enough OEIS fun; back to the task at hand.
While $G[1000,1000,1000,1000]$ contributes to $C[1001,1000,1000,1000]$, among other hyper-Catalans, so do other Geode elements, such as $G[1001,999,1000,1000]$.
Rubine, 2025~\cite{rubine2025hyper} explores this, taking the first step toward making the calculation of a single Geode element more efficient, with the following theorem.
\begin{theorem}[The 4D hyper-Catalan / Geode sum]
For a non-degenerate type vector $[m_2,m_3, m_4, m_5]  \ne [0,0,0,0]$ (at least one non-zero component),
\begin{align*}
C[m_2,m_3, m_4, m_5] &= G[m_2-1, m_3, m_4, m_5] + G[m_2,m_3-1, m_4, m_5]  \\ & \ \ 
+ G[m_2,m_3, m_4-1, m_5] + G[m_2,m_3, m_4, m_5-1]
\end{align*}
where $G[m_2,m_3,m_4, m_5]\equiv [ t_2^{m_2}  t_3^{m_3} t_4^{m_4} t_5^{m_5} ] \Gb$.
\end{theorem}
Clearly we have only natural powers in $\Gb$, so if $m_k=0$, then any $G$ indexed with $m_k-1$ will be zero, and may be omitted from the sum.
\begin{proof}
Theorem \ref{thm:geodefactorization} says $$
-1 + \Sb[ t_2,  t_3, t_4, t_5] = (t_2 + t_3 + t_4 +t_5) \Gb[t_2,  t_3, t_4, t_5].
$$
Extracting the coefficient of $\T^\m$ from both sides, we ignore the $-1$ as $\m$ is non-zero:
\begin{align}
\nonumber
[ t_2^{m_2}  t_3^{m_3} t_4^{m_4} t_5^{m_5} ] \Sb[ t_2,  t_3, t_3, t_4]& = [t_2^{m_2} t_3^{m_3} t_4^{m_4} t_5^{m_5} ] \sum_{2 \le j \le 5}  t_j \Gb[t_2,  t_3, t_4, t_5]
\\ \nonumber  &  =   
[t_2^{m_2-1} t_3^{m_3} t_4^{m_4} t_5^{m_5} ] \Gb
+[t_2^{m_2} t_3^{m_3-1} t_4^{m_4} t_5^{m_5} ] \Gb
\\[-2pt] \nonumber   & \ \ +[t_2^{m_2} t_3^{m_3} t_4^{m_4-1} t_5^{m_5} ] \Gb
+[t_2^{m_2} t_3^{m_3} t_4^{m_4} t_5^{m_5-1} ] \Gb
% \\ \nonumber
% C[m_2,m_3, m_4, m_5] &= G[m_2-1, m_3, m_4, m_5] + G[m_2,m_3-1, m_4, m_5] 
% \\ & \ \ + G[m_2,m_3, m_4-1, m_5] + G[m_2,m_3, m_4, m_5-1]
\end{align}
from which we conclude the theorem.
\end{proof}
Notice each $G$ has the same sum of array indices, and that the sum for $C$ is one more than that.  
We can solve for a Geode element. We add 1 to $m_2$:
\begin{theorem}[The 4D Geode Recurrence]\label{thm:4DGeode}
\begin{align*}
 G[m_2, m_3, m_4, m_5]  & = C[m_2+1,m_3, m_4, m_5] - G[m_2+1,m_3-1, m_4, m_5] 
 \\[-2pt] & \ \ - G[m_2+1,m_3, m_4-1, m_5] - G[m_2+1,m_3, m_4, m_5-1]
\end{align*}
\end{theorem}
Rubine showed that when $G[m_2,m_3,m_4,m_5]$ has only one non-zero $m_k$, then the element is a Fuss number~\cite{Fuss1795}, which is a particular kind of hyper-Catalan number that counts subdivisions into a single kind of shape (all quadrilaterals, for instance).  In particular, from Theorem \ref{thm:4DGeode} it is easy to see:
\begin{equation}
 G[m_2, 0,0,0]   = C[m_2+1, 0,0,0],
 \end{equation}
 a Catalan number.

Rubine goes on to describe a general recursive algorithm for computing Geode elements where at each step the index $k$ to be incremented has the largest $m_k$, and the first such in the event of a tie.
In the special case when we start from a diagonal element $G[n,n,n,n]$, that index will always be $2$, the first element, which we assume here.
For simplicity, we assume in the code samples that it is always the first index we will be incrementing (which is index 0 in the code; apologies for the mismatch).

The system used is largely freeware: Python, in particular Sympy~\cite{Sympy2017} bignums, in a Jupyter Notebook in a Chrome tab on an old PC shamefully running Windows 10 past its expiration date.

\begin{figure}
    \caption{Straightforward hyper-Catalan and Geode calculation}
    \label{fig:code1}
    \begin{lstlisting}[language=Python] %,frame=single]
    
from sympy import *
def C(m):
    return factorial(E(m)-1) / (Fact(m) * factorial(V(m)-1))
def V(m):
    return 2 + sum((i+1)*m[i] for i in range(len(m)))
def E(m):
    return 1 + sum((i+2)*m[i] for i in range(len(m)))
def Fact(m):
    return prod(factorial(m[i]) for i in range(len(m)))
def G(m):
    m[0] += 1
    s = C(m)
    for i in range(1,len(m)):
        if m[i] > 0:
            m[i] -= 1
            s += -G(m)
            m[i] += 1
    m[0] -= 1
    return s
    \end{lstlisting}
\end{figure}

The code in figure \ref{fig:code1} uses the general Geode recurrence to compute arbitrary Geode elements (not just 4D).
(In all the code samples, we eliminate error checking and log printing for clarity and brevity.)
This code is faster than the naive polynomial division, but it is still rather slow.

Writing the repeated subscripts is getting old; let's abbreviate the 4D diagonal Geode elements as $H(n) \equiv G[n,n,n,n]$.  Using Python's \%\%time, on my aging PC we determine that computing $H(1)=12344$ takes 0 seconds (i.e., is too quick  to measure), $H(2)=2408941884$ takes .004 seconds, $H(3)=894971463204720$ takes .060 seconds, $H(4)=446324644841317281200$ takes 1.26 seconds, $H(5)=$ $263656050352833337510832640$ takes 28.4 seconds, and $H(6)=\newline 173882340006327290808417397911384$ takes 654 seconds.
It is not that obvious what the complexity of the computation of $H(n)$ is, but it is obvious that this implementation of the Geode recurrence will not make it to $H(1000)$ anytime soon.

\section{Caching}

At those speeds we do not have much hope to compute $G[1000,1000,1000,1000]$ in a timely manner.
For our first speedup attempt, we add caches to save evaluations of $G$ and $C$.

There was substantial experimentation at this stage.
It turns out the hyper-Catalan cache is unnecessary; each $C$ element used is only used once in the calculation of $H(n)$.

We see each Geode calculation touches three Geode neighbors on a 4 dimensional lattice.  
So, except for edge cases, each calculated Geode element is used three times.
We experimented with deleting each Geode cache entry after it was no longer required, assuming memory was a bottleneck, but it did not seem to make any difference, at least up to $H(300)$.

We also experimented with using Sympy's \code{fallingfactorial} rather than the quotient of factorials in calculating $C$.
That was less effective than adding a factorial cache, and dividing cached values. 
The factorial cache was well-utilized; for $H(300)$ there were 164 million calls to factorial, but only 700 distinct inputs.

We can see that the Geode cache for one $n$ does not help for different $n$s.  
That is obvious from our observation about the sum of the indices $F_\m= m_2+m_3+m_4+m_5$; for the $G$s of interest in calculating $H(n)$, that sum is always $4n$.
Similarly, for the $C$s the sum is $4n+1$, meaning we are concerned with the set of subdigons with exactly $4n+1$ faces.
But, as we will see in the next section, not that entire set.

Figure \ref{fig:code2} shows the changes to add the caches; unchanged functions are not repeated.
The initial caching attempt was able to compute $H(200)$, 1201 digits, in 27 minutes, and $H(250)$, 1503 digits, in 20 hours.
The version shown completed $H(300)$, 1805 digits, in 5.9 hours, $H(400)$, 2409 digits, in 25.6 hours,  and $H(500)$, 3014 digits, in 82 hours.
\newcommand{\bigO}{\mathcal{O}}

Manuel Kauers (personal communication) suggests this implementation of $H(n)$ is of complexity $\bigO(n^5)$.
The number of arithmetic operations grows at $\bigO(n^3)$, but as the length of the numbers in the computation is proportional to $n$, so we get an additional factor of $n^2$ due to the bignum multiplication and division calculations.
Scaling up the $H(500)$ time by $2^5$ means we expect $H(1000)$ to take around 109 days with this code, not out of the question, but we want something faster.

\begin{figure}
    \caption{Calculation with cached Geode and factorials}
    \label{fig:code2}
    \begin{lstlisting}[language=Python] %,frame=single]
    
Gcache = { } ; Fcache = { } 
def MYfactorial(n):
    global Fcache
    if n not in Fcache:
        Fcache[n] = factorial(n)
    return Fcache[n]
def Fact(m):
    return prod(MYfactorial(m[i]) for i in range(len(m)))
def C(m):
    return MYfactorial(E(m)-1)/(MYfactorial(V(m)-1)*Fact(m))
def mstr(m):  # turn array indices into a dictionary key
    return ','.join([str(mi) for mi in m])
def G(m):
    global Ccache
    key = mstr(m)
    if key in Gcache:
        return Gcache[key]
    m[0] += 1
    s = C(m)
    for i in range(1,len(m)):
        if m[i] > 0:
            m[i] -= 1
            s += -G(m)
            m[i] += 1
    m[0] -= 1
    Gcache[key] = s
    return s
    \end{lstlisting}
\end{figure}

\section{Geode Slices for Computational Efficiency}

I thought to make another major speedup, I was either going to have to master the magical math in Amderberhan, \etal{}~\cites{Amdeberhan2025,Zeilberger2025} or I would have to understand what my recursive Geode calculation was doing in more detail.
I bet on the latter as the most expedient path to $G[1000,1000,1000,1000]$, but I am confident that the next challenge, $G[1000,1000,1000,1000,1000]$, is going to need the math.

Each Geode calculation uses three other Geode elements, neighbors in the 4D Geode lattice, and one hyper-Catalan element.
Figure \ref{fig:H2} shows the Geode elements involved in computing $G[2,2,2,2]$, a small case to be sure, but sufficiently illustrative. 
We see we get a $3 \times 3 \times 3$ lattice; in general the lattice for $H(n)$ will be $(n+1)\times(n+1)\times(n+1)$.
The figure follows the convention that trailing zeros of the type are omitted.

Each Geode element requires three Geode elements to compute, the one below it, the one to the left (and a bit down) and the one diagonally down.
The lattice confusingly has lines through nodes that are not necessarily connections, e.g. $G[3,2,1,2]$ is only connected to nodes more toward the origin $G[8]$, namely  $G[4,1,2,2],  G[4,2,0,2], $ and $ G[4,2,1,1]$, and not connected to $G[3,2,2,1]$.

We define a \textbf{slice} as one of the two dimensional slices of the lattice with a constant $m_3$; Slice $s$ has $m_3=s$.
The important observation is to compute Geode elements in slice $s$ only requires two `lower' Geode elements in slice $s$ and one Geode element in slice $s-1$.
Our plan is forming: compute one slice at a time, use it to incrementally compute the next slice.

\begin{figure}
    \centering
    \includegraphics[width=0.8\linewidth]{GeodeH3.png}
    \caption{The lattice of Geode elements used to compute $G[2,2,2,2]$; the number before the colon is the initial recursion depth.}
    \label{fig:H2}
\end{figure}

We'll divide the code into multiple parts; figure \ref{fig:code3a} shows the Geode computation logic.
The goal is to produce a Gcache, a Python dictionary of an entire slice's Geode values.
The input is the slice number $\code{s}$, the hyper-Catalan cache for that slice $\code{Ccache}$, and the Geode cache from slice $\code{s-1}$, called $\code{oldGcache}$, used when decrementing $m_3$.
The for loop over types $\m$ generates all the $\m$s in the current Geode slice in the order we need to calculate them without any recursion, bottom left to top right in each slice; we note the $m_k$ always sum to $4n$.
The code would fail with index errors if the order was incorrect.

The rest is the usual Geode Recurrence, getting the hyper-Catalan elements from the Ccache and the Geode dependencies from one of the two Gcache. 

\begin{figure}
    \caption{Geode Slice Calculation}
    \label{fig:code3a}
    \begin{lstlisting}[language=Python] %,frame=single
    
def MakeGcacheSlice(s, n, Ccache, oldGcache):
    Gcache = { }
    for m in [ [4*n - (s + m4 + m5), s, m4, m5] 
                for m4 in range(n+1) for m5 in range(n+1)]:
        m[0] += 1
        r = Ccache[mstr(m)]
        for i in range(1,len(m)):
            if m[i] > 0:
                m[i] -= 1
                k = mstr(m)
                r -= oldGcache[k] if i==1 else Gcache[k]
                m[i] += 1
        m[0] -= 1
        Gcache[mstr(m)] = r
    return Gcache
    \end{lstlisting}
\end{figure}

\section{Hyper-Catalan Ratios}

We need to efficiently calculate the hyper-Catalans as they contribute to the Geode elements.
Again we do it on a slice by slice basis.

\begin{figure}
    \centering
    \includegraphics[width=1\linewidth]{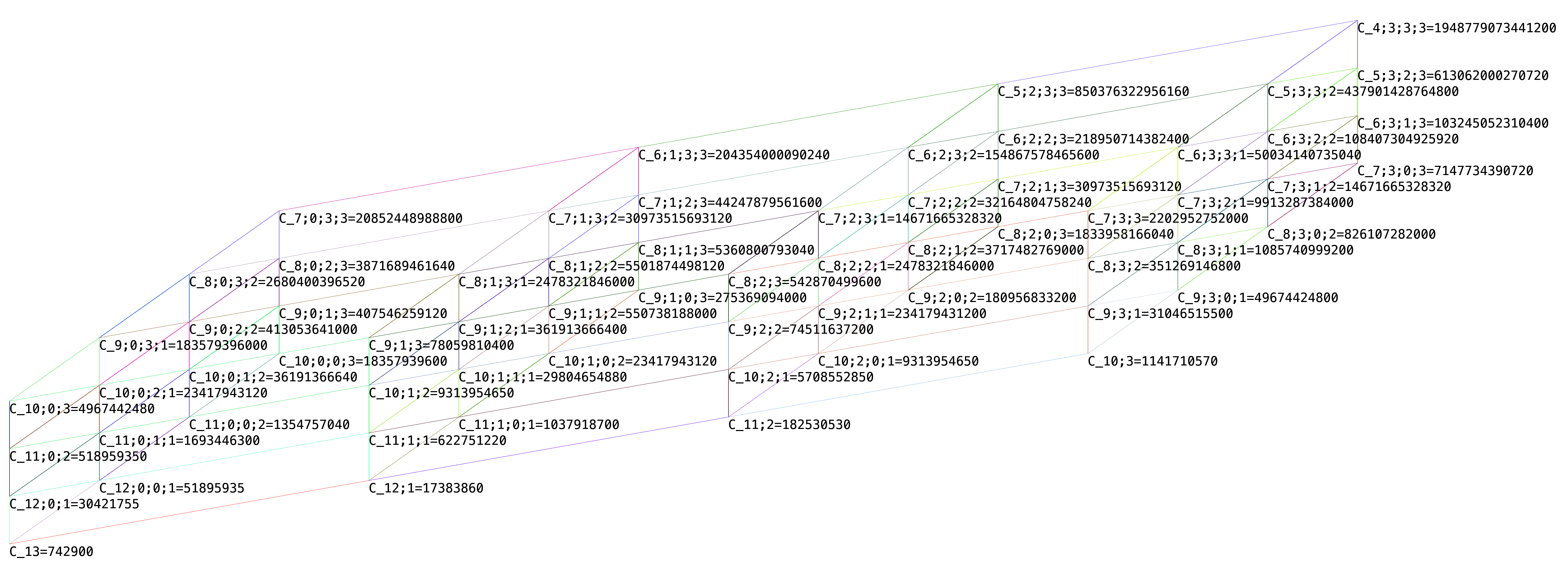}
    \caption{The hyper-Catalan Slice structure for $G[3,3,3,3]$}
    \label{fig:CH3}
\end{figure}

The hyper-Catalan lattice shown in figure \ref{fig:CH3} is the same shape as its corresponding Geode lattice.  The hyper-Catalan lattice indices add up to $4n+1$ because $m_2$ is one more than its respective Geode element.

We have two tasks ahead of us: Efficiently computing hyper-Catalan Slice 0, an $(n+1)\times(n+1)$ array, and efficiently computing hyper-Catalan Slice $s$ from hyper-Catalan Slice $s-1$.
Both rely on the same, more general code, but to avoid new notation here we write out the three cases we use separately:
\begin{theorem}[Ratio of hyper-Catalan neighbors]
\begin{align}
    \dfrac{C[m_2 -1, m_3 + 1, m_4, m_5]}{C[m_2,m_3,m_4, m_5]} &=  \dfrac{m_2 E_{\m} }{(m_3 + 1) V_\m}    \label{eqn:casem3}
\\ 
 \dfrac{C[m_2 -1, m_3 , m_4+1, m_5]}{C[m_2,m_3,m_4, m_5]} &=  \dfrac{m_2 E_{\m} (E_{\m} + 1) }{(m_4 + 1) V_\m (V_\m+1)}
\\ 
 \dfrac{C[m_2 -1, m_3 , m_4, m_5+1]}{C[m_2,m_3,m_4, m_5]} &=  \dfrac{m_2 E_{\m} (E_{\m} + 1)(E_\m+2) }{(m_5 + 1) V_\m (V_\m+1)(V_\m + 2) }    
\end{align}
\end{theorem}
%    t =  RisingFactorial(Em, k-j) * m[j-2] 
%    b = RisingFactorial(Vm, k-j) * (1 + m[k-2]) 
\begin{proof}
We do the middle one and leave the others to the reader.

Abbreviate the numerator type $\n= [m_2 -1, m_3 , m_4+1, m_5]$.
We have: 
\begin{align}
% E[m_2 -1, m_3 , m_4+1, m_5] &= 
E_\n &=  1+2(m_2-1) + 3 m_3 + 4(m_4+1) + 5m_5 = 2 + E_\m,
\\ 
% V[m_2 -1, m_3 , m_4+1, m_5] &=
V_\n &=  2+(m_2-1) + 2m_3 + 3(m_4+1) + 4m_5 = 2 + V_\m,
\\
\n! &= (m_2 -1)! \, m_3! \, (m_4+1)! \, m_5 = (m_4+1) \m! / m_2
% \end{align}
% so
% \begin{align}
\\
C_\n &= \dfrac{(E_\n - 1)!}{(V_\n - 1)! \, \n !}  = \dfrac{m_2  (E_\m + 1)! \, }{ (m_4+1)  (V_\m + 1)! \, \m! \,}
\end{align}
\begin{align}
 \dfrac{C_\n}{C_\m} 
 = \dfrac{ m_2 (E_\m + 1)! \,  }{(m_4+1)(V_\m + 1)! \, \m! \,  }
 \cdot \dfrac{(V_\m - 1)! \, \m !}{(E_\m - 1)!} = \dfrac{m_2 E_\m (E_\m + 1) }{(m_4+1) V_\m (V_\m+1)}
 \end{align}
\end{proof}
We capture the general result in the code in figure \ref{fig:code3b}; the subtractions of 2 adjust the indices given as parameters (starting at 2) to the internal representation starting at zero. 

\begin{figure}
    \caption{Ratio of hyper-Catalan neighbors}
    \label{fig:code3b}
    \begin{lstlisting}[language=Python] %,frame=single

 def Crat(m, j, k): # ratio of C(m-e_j+e_k)/Cm 
    return RisingFactorial(E(m), k-j) * m[j-2] / (
            RisingFactorial(V(m), k-j) * (1 + m[k-2]) )
    \end{lstlisting}
\end{figure}

The function allows us to calculate hyper-Catalans while avoiding factorials of large numbers; in practice $1 \le k-j \le 3$ so the rising factorials are at most two multiplications each.

\section{Efficiently Calculating Hyper-Catalan Slices}

Now that we have access to the ratios between nearby hyper-Catalan elements in the lattice, it is relatively straightforward to calculate the needed hyper-Catalans.
Figure \ref{fig:code3c} shows how the first hyper-Catalan slice, Slice 0, is calculated.  The if statement chooses the appropriate previously calculated neighbor for the most efficient calculation.

The next piece of the puzzle calculates the general hyper-Catalan slice from the previous one (figure \ref{fig:code3d}); for successive slices we're always incrementing $m_3$ for the most efficient calculation, as then the rising factorials only have a single factor so do no calculation at all.

\begin{figure}
    \caption{Creating hyper-Catalan Slice 0}
    \label{fig:code3c}
    \begin{lstlisting}[language=Python] %,frame=single

def MakeCcacheSlice0(n):
    cr0 = [ [1 + 4*n - (m4 + m5), 0, m4, m5] 
             for m4 in range(n+1) for m5 in range(n+1) ]
    firstm = cr0[0]
    Ccache = { }
    Ccache[mstr(firstm)] = C(firstm)
    for m in cr0[1:]:
        m2, m3, m4, m5 = m[0], m[1], m[2], m[3]
        if m4 != 0:
            om = [m2+1, m3, m4-1, m5]
            crat = Crat(om, 2, 4)
        else:
            om = [m2+1, m3, m4, m5-1]
            crat = Crat(om, 2, 5)
        Ccache[mstr(m)] = Ccache[mstr(om)] * crat
    return Ccache    
\end{lstlisting}
\end{figure}

\begin{figure}
    \caption{Creating hyper-Catalan Slice r}
    \label{fig:code3d}
    \begin{lstlisting}[language=Python] %,frame=single

def MakeCcacheSlice(s, n, oldCcache):
    if s==0:
        return MakeCcacheSlice0(n)
    Ccache = { }
    cr = [ [1 + 4*n - (s + m4 + m5), s, m4, m5] 
             for m4 in range(n+1) for m5 in range(n+1) ]
    for m in cr:
        om = [m[0]+1, m[1]-1, m[2], m[3]]
        crat = Crat(om, 2, 3)
        Ccache[mstr(m)] = oldCcache[mstr(om)] * crat
    return Ccache
\end{lstlisting}
\end{figure}

That is pretty much it; we put it all together in figure \ref{fig:code3e}.
\begin{figure}
    \caption{Efficient Calculation of Geode Diagonal Elements via Geode and hyper-Catalan slices}
    \label{fig:code3e}
    \begin{lstlisting}[language=Python] %,frame=single

def H(n):
    Ccache = { }
    Gcache = { }
    for s in range(n+1):
        Ccache = MakeCcacheSlice(s, n, Ccache)
        Gcache = MakeGcacheSlice(s, n, Ccache, Gcache)
    return Gcache[mstr([n,n,n,n])]
    \end{lstlisting}
\end{figure}

We could get away with a simpler, less efficient calculation for the initial slice, but for $H(1000)$ that would involve the computation of over a million hyper-Catalan numbers, which would take some time.
The cost was a more complicated \code{Crat}; for all the other slices we only need the first ratio, equation \eqref{eqn:casem3}.

\section{G[1000,1000,1000,1000]}

With the Slice-based implementation, $H(100)$ took 0.010 hours, $H(200)$ took 0.10 hours, $H(500)$ took 3.9 hours, $H(600)$ took 7.0 hours, and $H(1000)$ took 35.9 hours, coming in at 6303 digits.  We report the value as:

$ \  $

\noindent $G[1000,1000,1000,1000]= {}$

\noindent$\seqsplit{140604899259853103845676834712574810460853234736221730236256365412996835013135856006441319510456938285864287568734450135724848685182227852036329630817239601420905927701870187661110631936586189788184954823392652746425867144997401927129384104962331330281480173214764429464004355906854189800523031871022084126921807323325494407384669927928116121553257930625502629389085078102583931806230141348175588783571473308262539688853463304259592404017092301040534094688671172003439664116559795350569547714718469070504741491635114629029579382847965813461698227276988417019745638633875765491351113372461203167145626597243401369052257725928780008822810342230174508686655925179948869969600540601831352358507450326381143639410017753426708965040437560341599940781153997270789865972953566837161217111253410152588454187031672086381148004444503331286985075614635134579715605324549132491425209355090931770155348695020935214597326397535881566089783460534754476927934607636247602008734973136407399613814356919349096060423885861064127209733788611149392791155046352241002596206352678058868697116210485667865768727138158349842089850815429182049383715155452482314099422681356039083366667206371829374680400863994340185337024428321294671278260248273485784753132775857523106522429224818898503958765097426127156859278072896344144860138862083397784918024637751245622595436812008981396305654928700295787285386413524064091665639298525607144863906995970164681494242167551343138751539104013398914135519752693759815500039688686904364745136102594311961512437677655841847751273736238050932997917582162121171007909073706411467773836002411560555598818685202207432864988491277060400736440823794505744892190727319798715341893114570479040407125436582220522541754888900927283933318540829135440295031956055975464860389827610205203900640143533258476950741512684000111837262814266675133780203684343557073989855293017302710451533142068085372635709871360423862275010407434205481390367417715912403798288034215716231688609673455417950565607963300641584412033752281584180632654995947340565489741567872635297326684890936099010597627865318664734736272708345399840442529835698209813581313765734063147507646003708000064491745055809494580288706003572002050153455663976862763316766331487316845772064037336512285665013150851614274438657298230944952429390979292138298166828516761224455331310932418981924950663744703120913255566749907294781811708155999126804710197336531329033262697682256137102635788787585787287032393701657901662681426881705815579622651296216961882391915838723777243415809080896700338488033053462525162966444110331768637704549911446123994893842308252159361001750284477664678298574312630390055150154862428967682720927248354150919159249402583625805517240564840756165143710645927055256449816532331135164310410236440073956354540208910481068211596024059942337778492835210222091266792900201362803641542579386136972503389245386037382408661905323636669108559357640398780182182722121697323079727216791999241933139100854935416367512033229293190688718672205910273763957084007195766101139365544105361149356604607216470711927133730760059757754081290275944384321738357172790574363985460625928502357342025883966373523755661797839685529875114297940801246181575798107460249232559289067533130235330441607471570370666004496847057556976257149067353883004273784329669544947008750262805244973275950945669353510185285942778997061200661778555911920803818755673117377528010202480346335503987346440846818035907176467028504571770857983721800681865996092508643963452913949120073960181721453521181436423996919635288221719762597339537366005343157458285151228236253952331602015124044315066410643067948057777484668006089824229119496335606170374895811068250555264051255273075212177057223205732458897280228654899396597246078280898233970055239924278112016542251101758636396219497382222869892800456254233784966965069342231544137478604784793853481640461959596247455096563180699528630869218683318785581002911035031671947228466605574733129250185243351913746070349291657885361550859607353993471495353269698241520421315041958394520561048447421336417602179505041930247109888715008530389554135809579069273452970359297858004910902488028423764624011587073779168867942376967832468302002820726711102456454812577256662726554949047285522805250254483294805135204680005913860930128057162631339038161679463566865197694695621416101448684963072055453247069888965819153787162929350597302748153868937615165344438496566730249578084632791178468637634220562188853879714489106914316592676539780160759439049250772587202038012097412930217080842538361717079826741952082930939679165877102310337357348657199073614314055610136594258676650887940206191102246782916881635894091898531364106279830710134664321519861645244099442958538150461293341329153129378014002767460774359622384991970390570111252052611815315096125251549452087798688047265944186689690368282502861420583644556147894931907287540527231638423712496799199705688340981436556977566489099123199126951366539481469177858321130968241283803032163229530771037428580866018530042590921373580797009300674696864904720866416852283321993700355402475805999761808751715742023520912165867566640288564382733712758291448595176374776056070904929622953030352679119536687289148724909857062365249702360804547740204718522503531452763024182117324718121267287547520332772497984152283204843771526656614306791483252856532197957213156986482249273000448405169570954770697503355192452680992631860613203801621715228798467647111777767129334168930597904948132441233768057031586482713250596652011956905761019477980316024037377774231001542193678659024962565764659362106682528269046597134754936236844062127351939107622346889113843502847363333179942908854995106770521168352461094395151721997522049227689119961010362176818658777492219959692050037613226880245832364551309965600514899997178000516091485010182761452630164428968926234029818712282056248437874461452730606214368568223504589887206135247825003901217605680022306251925827689527344978580259542328964519136205365029414092102503833549476633414239544648341512407713629184000}$.

\section{Future Directions}

Computing large four dimensional Geode diagonal elements is quite a challenge.
By understanding the recurrence lattice as a layering of slices, we are able to make the  calculation of $G[1000,1000,1000,1000]$ tractable.

Note the calculation presented here is not the recurrence sought by Amderberhan, Kauers and Zeilberger, which starts from a relationship between nearby Geode elements, and ends with a relationship betweek $H(n)$ and some lags $H(n-k)$, each expected to have a quite complicated rational function factor.

Even the two dimensional case is far from exhausted.
The known closed form covers more than $G[m_2,m_3]$, it covers $G[0,0, \ldots, 0, m_k, m_{k+1}]$ as well, so presumably a generalization of their recurrence is possible here, to the two consecutive shape case.
The case when we have two dimensions that are not adjacent, non-consecutive shapes, as in $G[m_2, 0, m_4]$, has not been addressed at all.
There is some evidence that this may not have the hyper-geometric properties of the consecutive shape result, namely the relatively large primes for the first three values of $G[n,0,n]$.
There are similar unexplored three dimensional subarrays as well.

There is no particular reason to focus on the larger dimension diagonal.
The calculation given here may be adapted to calculate non-diagonal Geode elements; I would be interested in how the approach of Amderberhan \etal{} could calculate those.

Amderberhan \etal{} jumped right to the more difficult problem, skipping the 
recurrences for the hyper-Catalan array, which is more seemingly hyper-geometric than the Geode, so the challenge of the unbounded number of variables can be addressed without the other difficulties.

For the four dimensional Geode calculation, the slice structure makes the $\bigO{(n^3)}$ arithmetic operations apparent, and presumably we still have an additional factor of $n^2$ to handle the large big number arithmetic, so complexity $\bigO{(n^5)}$.
We are not seeing the expected $\bigO(n^2)$ big number factor in the runs; for $H(1000)$, slices 1 to 11 took 23:07 and slices 990 to 1000 took 22:04, actually faster on the much larger numbers.  There is probably something here that could be understood better.

The remaining challenge given is to compute the five dimensional diagonal element $G[1000,1000,1000,1000, 1000]$.
That adds another dimension to the lattice, so another factor of $n$ to the complexity, $\bigO{(n^6)}$, so with the current technology would take around 1000 times longer than the 4d calculation, over four years.
We will need a new idea to complete that calculation in a reasonable amount of time.

\section*{Acknowledgements}

The author thanks Doron Zeilberger, Manuel Kauers, Tewodros Amdeberhan and Norman Wildberger for their encouragement and assistance. 

\section*{Declarations}

The author has no competing interests to declare that are relevant to the content of this article.
No funds, grants, or other support was received.
No datasets were generated or analyzed during the current study.

\bibliography{HyperCatBib.bib}

\end{document}